\documentclass[journal]{IEEEtran}
\IEEEoverridecommandlockouts
\makeatletter

\long\def\@makefntext#1{\parindent 1em\noindent\hbox to 1.8em{\hss\@makefnmark}#1}
\makeatother
\usepackage{cite}
\usepackage{amsmath,amssymb,amsfonts}
\usepackage{graphicx}
\usepackage{textcomp}
\usepackage{xcolor}
\usepackage{xspace}
\usepackage{booktabs}
\usepackage{amsthm}
\usepackage{algorithm}
\usepackage{algpseudocode}
\usepackage{xurl}

\definecolor{lightgray}{gray}{0.9}
\definecolor{darkergray}{gray}{0.8}
\definecolor{forestgreen}{RGB}{34,139,34}

\usepackage{hyperref}
\hypersetup{
    colorlinks=true,
    citecolor=green,
    filecolor=black,
    linkcolor=red,
    urlcolor=blue
}
\usepackage[capitalize]{cleveref}
\makeatletter
\makeatletter
\usepackage{titlesec}
\titlespacing*{\section}{0pt}{*1}{*0.5}
\titlespacing*{\subsection}{0pt}{*0.8}{*0.4}

\IEEEoverridecommandlockouts

\begin{document}

\title{Fast Relax-and-Round Unit Commitment with Sub-hourly Mechanical and Ramp Constraints}

\author{
\IEEEauthorblockN{Shaked Regev\IEEEauthorrefmark{1}, 
\thanks{Notice: This manuscript has been authored by UT-Battelle, LLC, under contract DE-AC05-00OR22725 with the US Department of Energy (DOE). The US government retains and the publisher, by accepting the article for publication, acknowledges that the US government retains a nonexclusive, paid-up, irrevocable, worldwide license to publish or reproduce the published form of this manuscript, or allow others to do so, for US government purposes. DOE will provide public access to these results of federally sponsored research in accordance with the DOE Public Access Plan (\href{https://www.energy.gov/doe-public-access-plan}{https://www.energy.gov/doe-public-access-plan}). Research sponsored by the Laboratory Directed Research and Development Program of Oak Ridge National Laboratory, managed by UT-Battelle, LLC, for the US Department of Energy.}
Eve Tsybina\IEEEauthorrefmark{1}, 
Slaven Pele\v{s}\IEEEauthorrefmark{1}}

\IEEEauthorblockA{\IEEEauthorrefmark{1}Oak Ridge National Laboratory\\
1 Bethel Valley Road, Oak Ridge, Tennessee, USA\\
Email: regevs@ornl.gov}
}

\maketitle

\begin{abstract}
We propose a novel computational method for unit commitment (\textbf{UC}), which does not require linearized approximation and provides several orders of magnitude performance improvement over current state-of-the-art. The performance improvement is achieved by introducing a heuristic tailored for UC problems. The method can be implemented using existing continuous optimization solvers and adapted for different applications. We demonstrate value of the new method in examples of advanced UC analyses at the scale where use of current state-of-the-art tools is infeasible. We expect that the capability demonstrated in this paper will be critical to address emerging power systems challenges with more volatile large loads, such as data centers, and generation that is composed of larger number of smaller units, including significant behind-the-meter generation.
\end{abstract}


\section{Introduction}
\label{sec:introduction}

Recent changes in the electric power industry indicate a growth in load, a decrease in the average size of generating units, and an increasing complexity in forecasting. According to the U.S. Energy Information Administration (\textbf{EIA}), over the past 50 years total electricity production across all sectors increased by a factor of 2.3, from 1.87 thousand TWh in 1974 to 4.3 thousand TWh in 2024~\cite{EIA_MER_2025_11}. 

This growth is expected to intensify due to the increasing digital share of the US economy. For instance, 2024 S\&P Global Market Intelligence Datacenters and Energy Report indicates that 43.8 GW of datacenter capacity are planned for deployment between 2025 and 2030\cite{spglobal2024}. At nominal capacity, data centers could contribute up to 384 TWh of annual electricity consumption, almost a 10\% increase in load from a single load category.

Generator EIA data shows they are becoming smaller. For US units in operation in 2024, excluding Alaska and Hawaii, average nameplate capacity for a generator built in 1974 was 205.22 MW. For a generator built in 2024 it was 49.5 MW~\cite{eia8602024}. For 2408 generating units listed in ``proposed'', the average nameplate capacity was 80.21 MW for units dated 2014-2024 and 193.83 MW for units planned for 2025-2030. However, the actual average capacity is expected to be lower. First, historically, proposed units between 2014 and 2024 were approximately 1.5-4.5 times larger than those ultimately built. 

Furthermore, EIA statistics can be assumed to be biased towards larger generating units belonging to utilities. A large part of new generating units may be coming from behind-the-meter generation and may be significantly smaller~\cite{FERC2025}. A basic assessment of the reviewed data shows that if the electricity demand grew by a factor of nearly 2.5 and the average generating unit size decreased by a factor of 4, the total number of generating units in the system, including those not accounted for by EIA, could increase by a factor of 10, based on past trends, and this trend is likely to continue.

In addition to size, there remains a question of time related system complexity. According to EIA, the share of generation from all thermal sources, which are traditionally considered predictable, increased from 83.7\% in 1974 to 93.1\% in 2004, but declined to 78.9\% in 2024~\cite{EIA_MER_2025_11}. This decline, were it to continue, would make the generation stack less predictable, contributing to increased operational uncertainty and greater computational and modeling complexity. 

Intra-hour constraints are becoming significantly more important, as evidenced by the IEEE Task Force on Solving Large Scale Optimization Problems in Electricity Market and Power System Applications~\cite{chen2023}. However, higher granularity of unit commitment involves more frequent unit commitment calculation: 4 times more frequent for the 15 minute interval currently used by many researchers, and 12 time more frequent for the 5 minute interval, which is often used for modeling non thermal resources.

Meanwhile, unit commitment formulations and solvers remain limited. Recent literature reviews show system sizes for existing academic studies, which tend to use CPLEX or Gurobi on standard production-grade hardware. The order of magnitude of the problem is usually in $100$s of generators (not buses). For instance, \cite{chen2023} discusses one study with 158 units, \cite{falvo2022} examines studies with 20-130 units, \cite{montero2022} surveys studies with 5-1870 units, with most of references in the range of 10-200 units, \cite{hong2021} reviews systems with a stated system size of up to 100 units. A survey of large unit commitment problems of 5-2709 generators is discussed in \cite{vanackooij2018}, with most reviewed studies using 100-200 units (Note: the original paper referenced by \cite{vanackooij2018} uses 2704 generators \cite{fu2007}). The studies directly targeting systems with finer temporal granularity show very limited scale: \cite{parvania2016} discusses a system with 32 generators for 30-minute unit commitment, \cite{dwyer2015} uses 54 units, 15-minute steps, \cite{kazemi2016} uses varying extents of granularity up to 5-minute intervals on a system with 10 dispatchable generating units. 

Studies directed at improving computational efficiency use comparably small simulations and one hour resolution. For example, \cite{kim2019} uses a bundle method on a 130-generator system, \cite{yu2024}  uses temporal bundles and 289 generators, and \cite{saavedra2020uc} demonstrates security constrained unit commitment on 69 generators. There are several reasons for the use of small system sizes in unit commitment problems. First, system scaling happens along both spatial and temporal granularity. Running 1000 unit commitment scenarios for a 150-unit system remains substantially faster than running one scenario for a 1000-unit system. Similarly, running 72 time decoupled unit commitment scenarios is faster than running one time-coupled scenario for a three days outlook. Second, existing unit commitment algorithms rely on mixed integer programs (MIPs) \cite{achterberg2013}. In the worst case, solving MIPs is NP-hard and requires time exponential in the problem size, to enumerate all possibilities. For more information about the mathematical reasons for the poor scaling behavior of MIPs, see optimization-specific reviews such as \cite{koch2022, morrison2016, berthold2015, ridha2020}. 

In earlier work \cite{Regev2025}, we proposed a fast approach to unit commitment for small, fast-ramping generating units, which are expected to play an increasingly important role due to the rapid expansion of data center infrastructure. In this study, we extend this approach to include larger and slower generating units found in today’s conventional grids. To achieve this, we first review the basic UC formulation from \cite{Regev2025} and then augment it with supplementary constraints usually found in literature~\cite{chen2023,montero2022}.

Our approach only accounts for the actual costs of generation in the current period, but is biased towards turning more efficient units on. It includes intertemporal runtime constraints and, later, ramping constraints. From that perspective, the algorithm is partially temporally decoupled and cannot be directly compared to fully coupled UC algorithms. However, it still presents an advancement over earlier methods because (i) it scales sub-quadratically rather than exponentially and (ii) does not linearize the problem, making solutions more accurate. Each extended formulation is evaluated from two viewpoints: algorithmic runtime and the 
objective function. For the relaxed unit commitment formulation, we further analyze the stochastic characteristics of the objective function.

The study is organized as follows. \cref{sec:review} summarizes the original model. \cref{sec:runtime} introduces intraday unit constraints. 
\cref{sec:ramp} extends the model to include ramping constraints. 
\cref{sec:nonconst} introduces a different ramping formulation and compares its results to \cref{sec:ramp}.
\cref{sec:summary} summarizes the results and suggests further research directions.

\section{Review of Original Unit Commitment Problem}
\label{sec:review}
\subsection{Mathematical Formulation}

In our previous work~\cite{Regev2025}, the problem formulation was designed to address a very specific system need, i.e. to accommodate a potentially large number of small- to medium-sized flexible resources, potentially with highly convex supply curves. This resulted in the objective function that emphasized cost structure compared to equipment or grid constraints. Each unit $i$ was modeled with quadratic cost curves and allowed to commit up to a certain volume between $P_{\min,i}$ and $P_{\max,i}$ as opposed to committing the entire $P_{\max,i}$ of a unit.

To address possible demand uncertainties, we introduce the reserve margin, similar to the robust approach discussed in literature \cite{montero2022,dwyer2015,zheng2015}. This ensures that if our largest committed generator goes out, and demand is higher than predicted, we can still supply it. Here we choose a simple Gaussian noise model, and will build on this in future work. 

The respective optimization problem for $n$ generators is:
\begin{subequations}\label{eq:OldRCUC}
\begin{align}
\min \quad &
\sum_{i=1}^n \left(a_i P_i^2 + b_i P_i + c_i\right)u_i
\label{eq:OldObj}
\\[6pt]
\text{s.t.} \quad & P_{\min,i} \le P_i \le P_{\max,i}, \; u_i\in \{0,1\}
\label{eq:varConstOld}
\\&
 \sum_{i=1}^n u_i P_i \ge D,
\label{eq:DemandOld}
\\
&
 \sum_{i=1}^n u_i P_{\max,i}
\ge D + 3\sigma_D + R
\label{eq:MaxConstOld}
\end{align}
\end{subequations}
The parameters supplied to the optimization solver are:
\begin{itemize}
    \item $a_i, b_i, c_i$: quadratic cost coefficients for each generator $i$.
    \item $P_{\max,i}, P_{\min,i}$: maximum and minimum power output for generator $i$. 
    \item $D$: current predicted demand.
    \item $\sigma_D$: standard deviation of demand forecast 
    \item $R = \max\{P_{\max,i}\}$: reserve margin.
\end{itemize}
The decision variables are:
\begin{itemize}
    \item $P_i \in [P_{\min,i}, P_{\max,i}]$ for $1\leq i \leq n$: power output of each generator $i$.
    \item $u_i \in \{0, 1\}$ for $1\leq i \leq n$: commitment status of generator $i$.
\end{itemize}

The constraint \cref{eq:DemandOld} ensures that we supply the projected demand. The constraint \cref{eq:MaxConstOld} ensures that if our largest generator goes offline  we will be able to supply demand if it is within $3\sigma_D$ of our prediction. This ensures we can still supply demand if any one generator goes out and our predicted demand is higher than expected. More complicated uncertainty models can fit into this framework, but developing them is beyond the scope of this work. 

We solve this problem by first relaxing $u_i \in \{0,1\}$ to $y_i\in[0,1]$ and solving the relaxed continuous problem. We reorder terms in the descending order in $y_j$ so that $y_j>y_{j+1}$. Then, we find the first $m$ generators such that $\sum_{j=1}^m P_{\max,j} \geq D + 3\sigma_D + \max_{j}[P_{\max,j}]$, and solve economic dispatch with all $m\leq k \leq n$, for a total of $n-m+1$. The solution is the subproblem that produced the lowest objective.

\subsection{Simulation and Results}
\begin{figure}[htbp]
\centering
\includegraphics[width=0.95\linewidth]{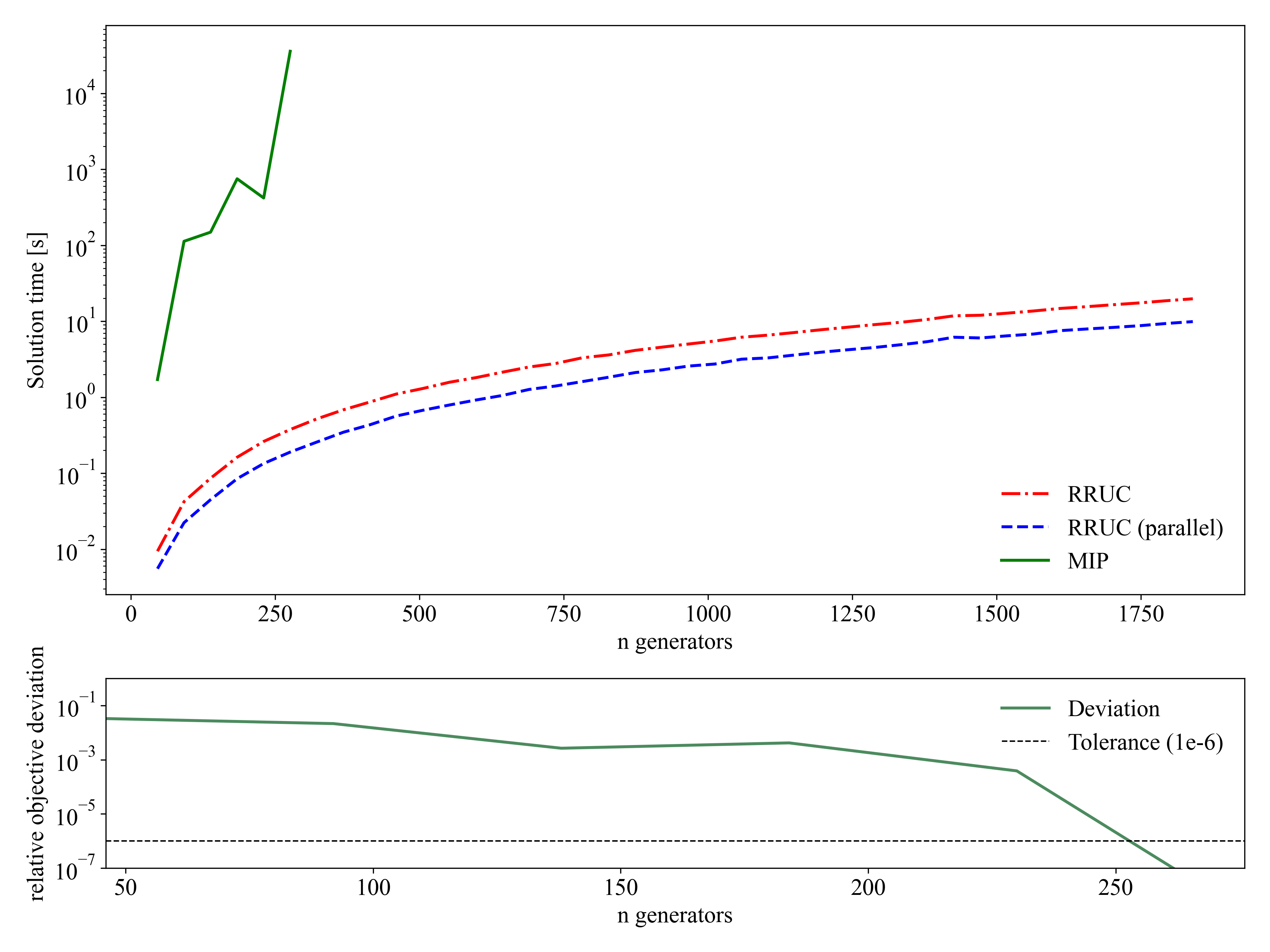}
\caption{\label{fig:old_time_obj}Test results for 46-1840 generating units \cite{Regev2025}. (Top) convergence time for both solvers, RRUC outperforms MIP by orders of magnitude. (Bottom) deviation in objective value calculated as (RRUC-MIP)/MIP. For larger problems, the deviation is below the tolerance set for the solvers.}
\end{figure}

To verify the performance of the approach, we used an increasingly large single node system with varied generators obtained from PJM Data Miner
\cite{pjm}. We fitted the piece-wise linear bids from PJM to quadratic cost curves with least squares, in line with \cite{yang2022,abdi2021}. 

We implemented our method in Julia~\cite{Bezanson2017Julia} with the MadNLP optimization package~\cite{Pacaud2024ExaModels}, which is based on our prior work~\cite{Regev2022HyKKT, Swirydowicz2025GPU}. As a reference MIP implementation, we choose the open-source Juniper package~\cite{Kroger2018Juniper}, which runs on the same JuMP platform~\cite{Dunning2017JuMP} as MadNLP. We ran both solvers on an Apple M3 Pro Laptop with 18GB RAM.

We highlight the key results of~\cite{Regev2025} in \cref{fig:old_time_obj}. Relax-and-Round Unit Commitment (\textbf{RRUC}) was $1,000$s of times faster than the MIP solver reference. We abandoned the latter after 270 units as it was taking about 10 hours on the available computational resources. 

The accuracy of the RRUC algorithm compared to MIP depended on the system size. On very small problems, a traditional MIP solver committed different units and achieved a smaller cost (objective function). The difference was under $6\%$. However, the difference decreased with the increasing problem size. On the 270-generator problem, the largest problem the MIP solver could solve, the difference in objective function was smaller than the tolerance set by the solvers and the two solvers selected the same units. 

\section{Unit Commitment with Runtime Constraints}
\label{sec:runtime}
\subsection{Mathematical Formulation}

We now expand the original formulation to include a wider variety of generators, some of which are slower large units found in utility power plants. Such units have runtime constraints and a maximum number of daily starts. Their startup and shutdown costs are a significant  portion of their total operating costs. The different components of runtime and other mechanical constraints are summarized in ~\cite{vanackooij2018,yang2022,tahanan2015}. 

We follow the general consensus using startup and shutdown costs and minimum on/off time constraints. Specifically, for each generator we add a minimum runtime $T_{\min,i}$ and a preprocessing stage before optimization. The preprocessing determines how long a generator has been running. If a generator is on but has run fewer periods than $T_{\min,i}$, the generator cannot be switched off. Such generators are considered must run and denoted as group $\mathcal{G}_m$. 

Further, each generator has a maximum number of startups during the 24-hour window, $S_{\max,i}$. The generators that are off and cannot yet be switched on due to the maximum number of starts are excluded from the optimization problem for that time period. The remaining generators are discretionary. We denote this group $\mathcal{G}_d$. The resulting optimization problem is:
\begin{subequations}\label{eq:OptimizationSimple}
\begin{align}
\label{eq:OptObj}
\min \quad &
\sum_{i \in G_m} \left(a_i P_i^2 + b_i P_i + c_i\right)+ \\
& \sum_{i \in G_d} \left[\left(a_i P_i^2 + b_i P_i + c_i\right) u_i
+ K_i(u_i - u_{t-1,i})^2\right]\nonumber
\\[6pt]
\text{s.t.} \quad & P_{\min,i} \le P_i \le P_{\max,i}, \; u_i\in \{0,1\}
\label{eq:varConst}
\\&
\sum_{i \in G_m} P_i + \sum_{i \in G_d} u_i P_i \ge D,
\label{eq:DemandConst}
\\
&
\sum_{i \in G_m} P_{\max,i} + \sum_{i \in G_d} u_i P_{\max,i}
\ge {D}_{\max,72} + 3\sigma_D + R,
\label{eq:MaxConst}
\\
&
\sum_{i \in G_m} P_{\min,i} + \sum_{i \in G_d} u_i P_{\min,i}
\le {D}_{\min,72} - \sigma_D\footnotemark .
\label{eq:MinConst}
\end{align}
\end{subequations}

\footnotetext{At the lower end, we plan for a smaller contingency ($\sigma_D$) than on the higher end ($3\sigma_D$). This is meant to reflect the asymmetric nature of load forecasts. It is far more likely that load will be much higher than predicted (corresponding to some large unexpected load), compared to much lower (corresponding to base load being removed). This also helps provide additional buffers in the event that more than 1 generator goes out, or uncontrollable generation is much lower than expected.}

The parameters supplied to the optimization solver are:
\begin{itemize}
    \item $a_i, b_i, c_i$: quadratic cost coefficients usually resulting from fuel costs or strategic bidding behavior.
    \item $P_{\max,i}, P_{\min,i}$: maximum capacity and minimum stable output, determined by a unit's physical characteristics
    \item $D$, $\sigma_D$: demand volume and the possible demand deviations in the robust model formulation.
 \item $K_i$: commitment change penalty for generator $i$\footnote{The true startup and shutdown costs are rarely the same. Averaging them prevents market distortion and removes artificial incentives to turn on or off.}.
    \item $u_{t-1,i}\in \{0, 1\}$: previous commitment status of unit $i$
    \item ${D}_{\min,72}$ and ${D}_{\max,72}$: minimum and maximum volumes of predicted demand in next 72 hours.
    \item $\mathcal{G}_m$: the set of must run generators.
    \item $\mathcal{G}_d$: the set of discretionary generators.
    \item $R = \max\{\max_{i \in \mathcal{G}_m} P_{\max,i}, \max_{i \in \mathcal{G}_d} P_{\max,i}\}$: reserve margin (changed to only include generators that can run).
\end{itemize}

The decision variables are:
\begin{itemize}
    \item $P_i \in [P_{\min,i}, P_{\max,i}]$ for $i \in \mathcal{G}_m$: power output of must-run generator $i$.
    \item $P_i \in [P_{\min,i}, P_{\max,i}]$ for $i \in \mathcal{G}_d$: power output of discretionary generator $i$.
    \item $u_i \in \{0, 1\}$ for $i \in \mathcal{G}_d$: commitment status of discretionary generator $i$.
\end{itemize}

Similar to the problem \cref{sec:review}, the problem defined in \cref{eq:OptimizationSimple} is a MIP, an NP-hard type problem. This problem is intractable and usually solved by a branch and bound method \cite{morrison2016}. We relax $u_i \in \{0, 1\}$ to $y_i \in [0, 1]$, reorder terms in the descending order, and determine which units are preferred to run. Since in our problem the $y_i$ variables represent which generators are committed and which are not, we find the first $\tilde{m}$ generators such that
\begin{equation}
\sum_{i\in\mathcal{G}_m} P_{\max,i}
+
\sum_{\substack{1\leqslant i \leqslant \tilde m\\ i\in\mathcal{G}_d}} P_{\max,i}
\;\ge\;
D_{\max,72} + 3\sigma_D + R_{\tilde m},
\label{eq:findmin}
\end{equation}
where
\[
R_{\tilde m}
=\max\!\left\{
\max_{i\in\mathcal{G}_m} P_{\max,i},
\;
\max_{\substack{1\leqslant i\leqslant \tilde m\\ i\in\mathcal{G}_d}} P_{\max,i}
\right\}.
\]
As  discussed above, some  previously committed units must operate due to runtime constraints. Therefore, $m\doteq |\mathcal{G}_m|+\tilde m$. We can calculate $m_{\max}$, the maximum number of generators that can operate, using \cref{eq:MinConst} analogously to the way we used \cref{eq:MaxConst} to derive \cref{eq:findmin}. This gives us \cref{eq:findmax}:
\begin{equation}
\sum_{i\in\mathcal{G}_m} P_{\min,i}
+
\sum_{\substack{1\leqslant i \leqslant \tilde{m}_{\max}\\ i\in\mathcal{G}_d}} P_{\min,i}
\;\ge\;
D_{\min,72} - \sigma_D,
\label{eq:findmax}
\end{equation}
with $m_{\max}\doteq |\mathcal{G}_m|+\tilde m_{\max}$. 
We then solve economic dispatch with all options of generators that can supply the load. Meaning, $\forall k: \; m \le k \le m_{\max}$, we set $y_{j\le k}=1$ and $y_{j>k}=0$ and solve the following problem:
\begin{subequations}\label{eq:OptimizationEco}
\begin{align}
\min \quad &
\sum_{j=1}^{k} \left(a_j P_j^{2} + b_j P_j + c_j+ K_j(u_j - u_{t-1,j})^2 \right)
\label{eq:EcoObj}
\\[6pt]
\text{s.t.} \quad &
P_{\min,j} \le P_j \le P_{\max,j}
\label{eq:EcoPower}
\\
&
\sum_{j=1}^{k} P_j \ge D.
\label{eq:EcoDemand}
\end{align}
\end{subequations}
Eqs. \eqref{eq:EcoObj}, \eqref{eq:EcoPower}, and \eqref{eq:EcoDemand} are the analogues of equations \eqref{eq:OptObj}, \eqref{eq:varConst}, and \eqref{eq:DemandConst} with the binary variables $u_i$ set to $1$ or removed, because we have already decided which generators to operate. Eqs. \eqref{eq:MaxConst} and \eqref{eq:MinConst} don't have analogues, because they only depend on $u_i$ (not $P_i$) and we selected $k$ so that they hold. These $m_{max} -  m + 1$ problems are completely independent and can be solved in parallel, though \cref{eq:OptimizationEco} solves much faster than \cref{eq:OptimizationSimple}. Then we take the solution which had the lowest objective over all values of $m \le k \le m_{\max}$.

\cref{alg:runtime} summarizes the workflow for this section.
\begin{algorithm}[htbp]
\caption{Runtime constrained RRUC.}\label{alg:runtime}
\begin{algorithmic}[1]
\For{$i=1$ \textbf{to} $\mathrm{number\_of\_periods}$}
  \State Update rolling 24-hour generator states (use them in steps 3 and 4).
  \State {Find must run units (with minimum on time not yet reached or with $P_{t-1,i}-P_{min,i}$ too large to shut off)}
  \State Find can't start generators (those that have exceeded  maximum allowed daily starts or are already ramping)
  \State Solve \cref{eq:OptimizationSimple} with $u_i\in\{0,1\}$ relaxed to $0\le y_i\le 1$.
    \State Order units in $\mathcal{G}_d$ into sets descending by $y_i$
  \For{each $\mathcal{G}_d$ that meets (\cref{eq:varConst})-(\cref{eq:MinConst})}
    \State Solve economic dispatch \cref{eq:OptimizationEco} with $\mathcal{G}_m, \mathcal{G}_d$ .
    \If{objective is lowest so far}
      \State Update best solution and objective.
    \EndIf
  \EndFor
  \State \textbf{save} Best solution and best objective.
\EndFor
\end{algorithmic}
\end{algorithm}


\cref{eq:OptimizationSimple} indicates that the problem has a weakly time-coupled objective with strongly time-coupled constraints. This approach has a potential downside of selecting a different set of generating units compared to a fully coupled approach, because the future payoffs and penalties are not accounted for at a given simulation step. It has the advantage of simplifying the computation by reducing the number of required variables. 

Results from \cref{alg:runtime} are not be directly comparable to fully coupled approach on computational time. Since we have shown our algorithm scales sub-quadratically~\cite{Regev2025}, we believe it is possible to proceed with this approach and expand it in the future to include a fully time-coupled objective. In this work, the runtime constraints are still time-coupled and the system is forced to converge while fully accounting for minimum runtime or the number of starts. 

Unlike~\cite{Regev2025}, we do not compare RRUC results to MIP results. The large system sizes, combined with new constraints, make MIP prohibitively computationally intensive for the available hardware (Apple M3 Pro, 18GB RAM). Instead we focus specifically on RRUC and adopt a new metric to assess the efficiency of the algorithm: percentage increase in the objective function and percentage increase in runtime as a function of increase in the system size. If the objective function and runtime do not explode as the system size increases, the algorithm maintains its scaling properties and is efficient.

\subsection{Simulation and Results}

The system was simulated as follows. The historical PJM peak was recorded on June 23, 2025 at 5 pm, with total load of 160.2 GW ~\cite{pjm}. We apply this value to scale the system using a 20\% cold reserve margin for resource adequacy, resulting in the system size of approximately 192 GW. Our selection of generating units is designed to represent the typical PJM generators plus the additional small-scale generators which we expect to be prevalent in the future as we explain in \cref{sec:introduction}. It includes 42 real PJM generators and amounts to 9047.9 MW  capacity for the entire system, which we will refer to as $P_{\max}\doteq \sum_i P_{\max,i}$. 

We multiply this generator set by 22 to pass 192 GW, using a random seed to rescale each generator in order to avoid the symmetry effects~\cite{montero2022}. The noisy copies of the original 42 generators have a shift in capacity parameters of $\pm 10\%$ and a shift in daily starts or runtime hour parameters of $\pm 1$. For various system sizes, we further use multiples of this system.

We solve the runtime constrained problem over 2304 time periods at 5 minute intervals (8 days). This covers the peak and adjacent days of 6/17/2025-6/24/2025. The latter four days of this period was the interval in which PJM experienced the highest change of load in 2025, with the night low of 78.6 GW and the day peak of 160.2 GW, with intraday changes of 55-64 GW. We start the simulation four days before the peak, so that our 72 hour forecasts have time to impact our unit commitments before we hit the peak. The first day's results are used as a warm-up for the system (so that it is not starting with all units off) and we exclude them from our cost function. We assume that the standard deviation of the demand forecast $\sigma_D=1\text{GW}$ for the full PJM system (with $924$ generators) and scaled accordingly for larger systems.

To incorporate the runtime constraints, we augmented PJM generator data with runtime parameters summarized in \cref{tab:datarun}. PJM data had no clear assignment of units to gas or coal fuel, so we manually classified the units using the following rules: No unit larger than 1000 MW was considered gas because, according to EIA 860~\cite{eia8602024}, there were no operable gas units of this size in PJM territory. Coal units were identified as having a minimum runtime of 24 hours, unless they had the starting price in the supply curve above USD 10/MWh. The cost-based fuel assignment relies on the assumption that coal units have lower operating costs than gas units, although coal prices can be above gas prices on equivalent energy content~\cite{eiam}. 

Cold, warm, and hot start costs are available from PJM~\cite{pjm}. Those may be subject to strategic bidding, so it was helpful to validate them against a different system, such as the EU~\cite{diw2013}. By normalizing warm and cold start costs by hot start cost, we obtain relative costs for coal and gas units that are broadly consistent: cold start is typically 1.4–2.3× the cost of hot start, and warm start is 1.1–1.6× the cost of hot start.

\begin{table}[htbp]
\caption{\label{tab:datarun} Runtime Assumptions and Data Sources}
\begin{tabular}{lrrr} \toprule
   Parameter & Source & Coal & Gas
\\ \midrule
Min runtime [min] & ~\cite{pjm} & 240-1440 & 0-1440
\\ Max daily starts & ~\cite{pjm} & 1-3 & 1-24
\\ Hot start [min] & ~\cite{diw2013,GE} & 60-240 & 25-120
\\ Warm start [min] & ~\cite{diw2013} & 120-480 & 60-240
\\ Cold start [min] & ~\cite{diw2013} & 360-720 & 120-300
\\ Relative cost of start, warm/hot & ~\cite{pjm,diw2013} & 1.19-1.42 & 1.12-1.58
\\ Relative cost of start, cold/hot & ~\cite{pjm,diw2013} & 1.74-1.93 & 1.35-2.25
\\ \bottomrule
\end{tabular}
\end{table}

\begin{figure}[htbp]
\centering
\includegraphics[width=0.95\linewidth]{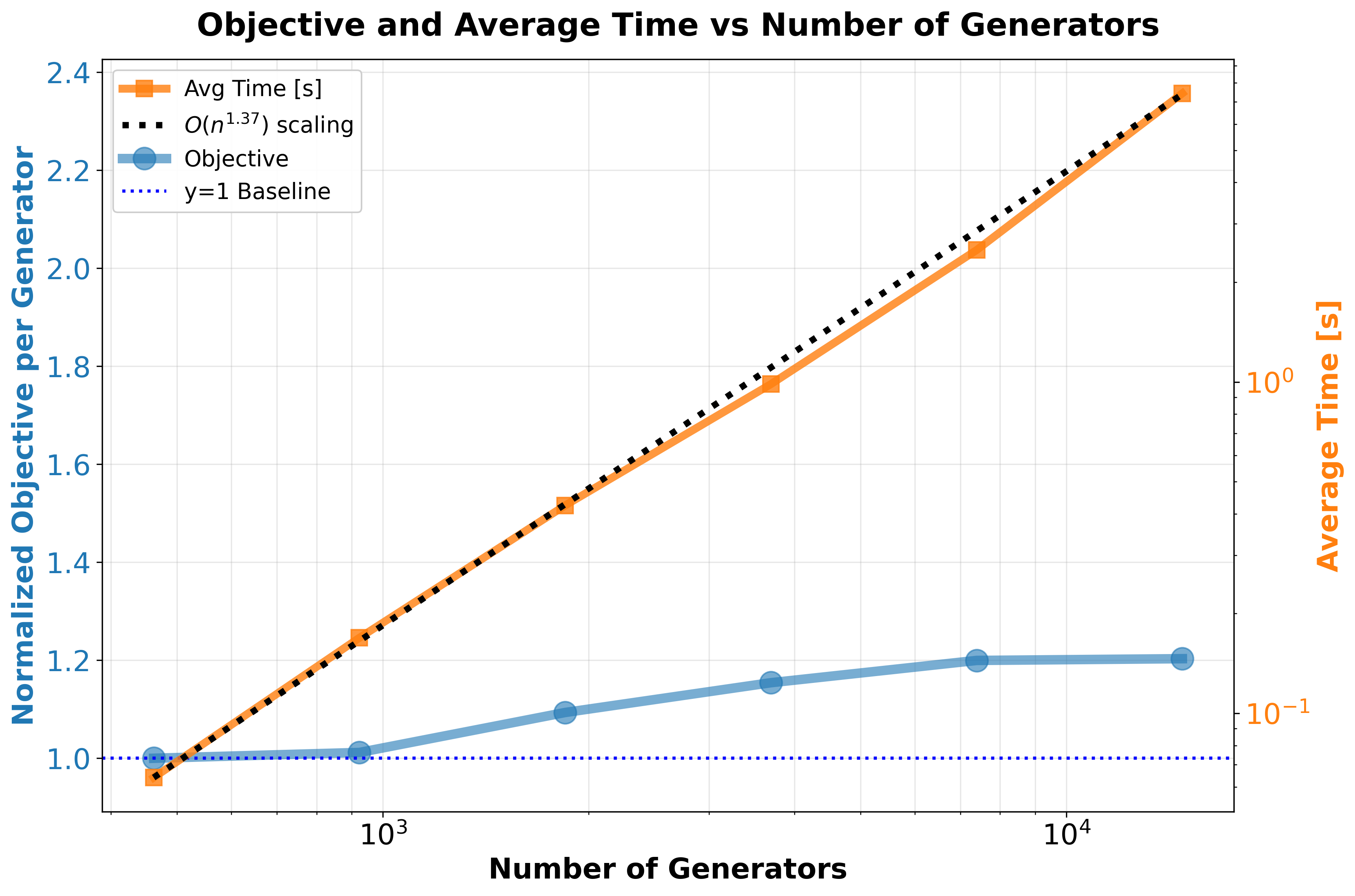}
\caption{\label{fig:mech_time_obj} Runtime and objective function value for RRUC algorithm including runtime constraints with 462-14784 generators. With $n$ as the number of generators, the runtime scales as  $\approx O(n^{1.37})$.  The objective function per generator increases less than 5\% for each doubling of the system size.}
\end{figure}

We present the results of the simulation in \cref{fig:mech_time_obj} using the same software and hardware as before. Since the reference MIP solver is too slow to be used as a benchmark for RRUC, in this study we are mainly interested in analyzing the scaling and performance of RRUC. The runtime and the objective function increase slightly faster than linearly with the number of generators. Specifically, the runtime increases on average about $2.6\times$ for every doubling of the system size. The $\approx O(n^{1.37})$ shows that the resulting systems are less complicated than those
arising from a typical optimization problem which have $O(n^{1.5})$ scaling. For a general optimization problem, $O(n^{1.5})$ scaling is believed to be ideal, due to the linear solver bottleneck from the ~\cite{morrison2016}.  

For the objective function, we rely on the per unit measure, or the average total cost, because our demand is scaled accordingly. The per-unit objective function increases less than 5\% for every doubling of the system. It is possible we could avoid the small loss of solution quality by increasing the run time.  An ideal solver's objective per unit would not increase with the system size. However, this model is just a stepping stone to the more realistic one in~\cref{sec:ramp}, so we will not spend time optimizing it. This indicates that larger systems do not result in significant diseconomies of scale and the algorithm remains cost efficient.

\section{Unit Commitment for Ramp and Sub Hourly Intervals}
\label{sec:ramp}
\subsection{Mathematical Formulation}


We now extend the problem to larger and more varied thermal power systems, which have have non-thermal units with intermittent characteristics. To operate such systems, we must account for system behavior at a higher resolution, which requires a sub hour granularity of unit commitment. Ramping is another important aspect of generator behavior at high temporal resolutions. 

There is a consensus among unit commitment studies that ramp is important in several ways \cite{chen2023,parvania2016,dwyer2015}. First, it presents a constraint on how quickly a unit may be expected to enter operation. Second, a unit that has been connected to the system but has not reached $P_{\min,i}$ still produces a certain volume of active power. That active power enters the grid and changes the net demand that is left for other units. Accurately accounting for ramp output results in different objective function values and different sets of committed units. 

\begin{figure}[htbp]
\centering
\includegraphics[width=0.95\linewidth]{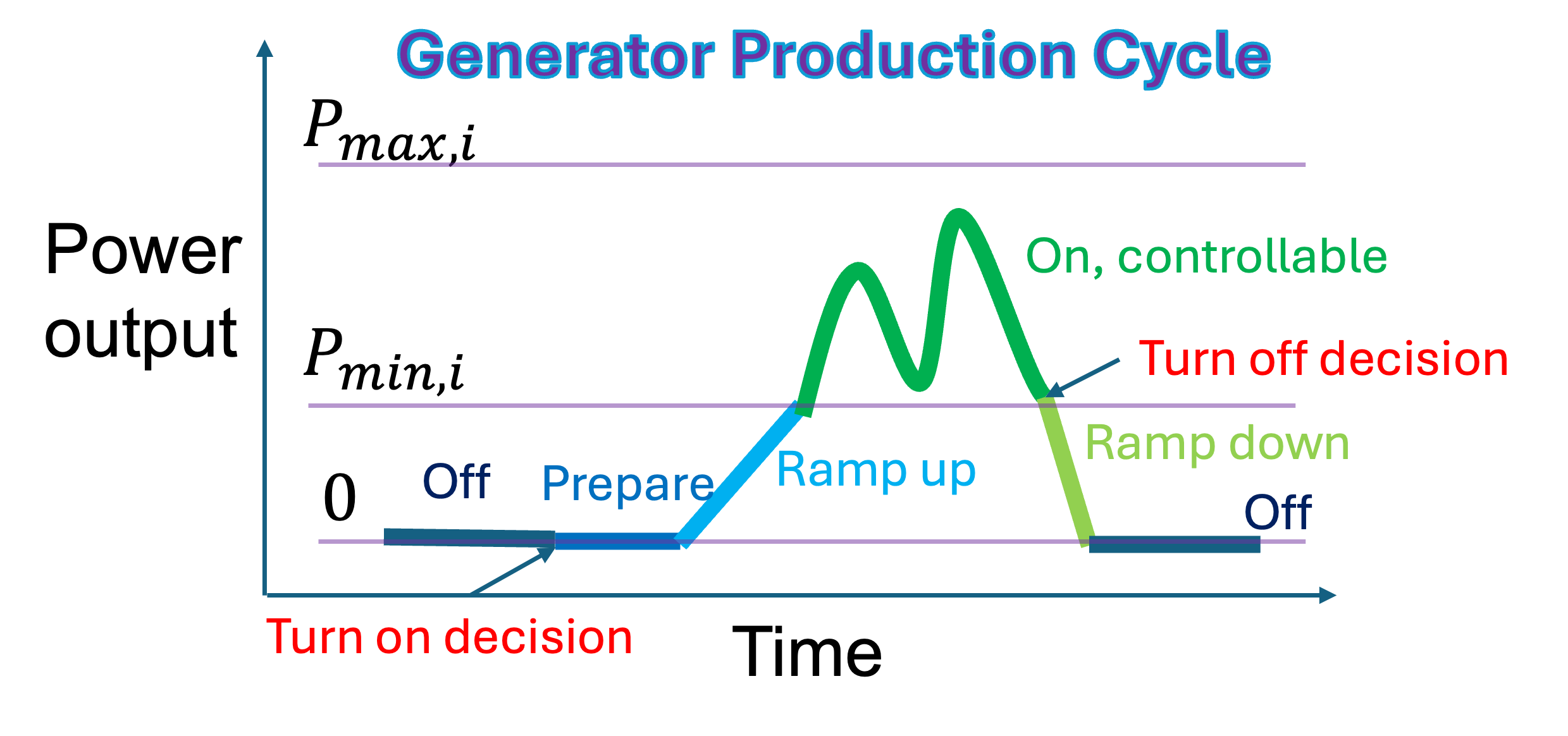}
\caption{\label{fig:ramp_profile} A typical production cycle for generator $i$ in \cref{sec:ramp}. When the generator is on, it can change its production power  between $P_{\min,i}$ and $P_{\max,i}$, provided it does not exceed its ramping power. It can also turn off if it is close enough to $P_{\min,i}$ compared to its ability to ramp down. The only other decision can be made when the generator is off, to prepare for ramping up. During all other periods, the generator is uncontrollable.}
\end{figure}

To address the ramp constraints in the model, we use a different set of variables. Instead of the original $u_i$ that indicated that a generator is on or off, we allow the generator to be off, preparing (for ramp up), ramping up, on, or ramping down. For the model purposes, we assume that a generator can only change its status in one direction. If a generator is ramping up or preparing, it is assigned a variable $v_i=1$. After the required ramp time in each stage, its status will switch to on. If a generator is ramping down, it will eventually switch off and enter the set of generators that can start. 

\cref{fig:ramp_profile} shows a typical production cycle of a generator and the allowed transitions. Each generator can only be in one state in a given time period, so we round ramping times from the data to whole periods. We round up to ensure that our solution is for a problem more restrictive than the actual problem, which means it is necessarily feasible. This means a generator must occupy each state for at least one time period. Similarly to \cref{sec:runtime}, we use preprocessing to account for the accumulated runtime and ramping periods. The information about the duration of on, off, or ramping is stored for each generator and used as an optimization parameter. 

The mathematical model used for ramp constrained sub hourly unit commitment is presented in \cref{eq:OptimizationEq5}.
\begin{subequations}\label{eq:OptimizationEq5}
\begin{align}
\min \quad &
\sum_{i \in G_m} \left(a_i P_i^2 + b_i P_i + c_i\right) \nonumber\\
&+ \sum_{i \in G_d} \left[\left(a_i P_i^2 + b_i P_i + c_i\right) u_i
+ K_i(1-u_i+v_i)^2\right]
\nonumber\\
& + \beta\sum_{i \in G_d} v_i \left(a_i P_{typ,i} + b_i + c_i/P_{typ,i}\right) 
\label{eq:OptObj5}
\\[6pt]
\text{s.t.} \quad & \max(P_{\min,i},P_{t-1,i}-r_{d,i}) \le P_i   \nonumber\\ &\le\min(P_{\max,i},P_{t-1,i}+r_{u,i})
\label{eq:cont5}
\\&
u_i\in \{0,1\}, \; v_i\in \{0,1\}, \; u_i+v_i\leq 1
\label{eq:bools5}
\\&
\sum_{i \in G_m} P_i + \sum_{i \in G_d} \left[u_i P_i + (1-u_i)(P_{\min,i}-r_{d,i})\right]  \nonumber
\\&\ge D-S_r,
\label{eq:DemandConst5}
\\
&
\sum_{i \in G_m} P_{\max,i} + \sum_{i \in G_d} (u_i+v_i) P_{\max,i} \nonumber\\
&\ge {D}_{\max,72} + 3\sigma_D - S_{\max,r},
\label{eq:MaxConst5}
\\
&
\sum_{i \in G_m} P_{\min,i} + \sum_{i \in G_d} (u_i+v_i) P_{\min,i} \nonumber\\
&\le {D}_{\min,72}- \sigma_D - S_{\min,r}.
\label{eq:MinConst5}
\end{align}
\end{subequations}

The decision variables are:
\begin{itemize}
    \item $P_i \in [P_{\min,i}, P_{\max,i}]$ for $i \in \mathcal{G}_m$: power output of must-run generator $i$
    \item $P_i \in [P_{\min,i}, P_{\max,i}]$ for $i \in \mathcal{G}_d$: power output of discretionary generator $i$
        \item $u_i \in \{0, 1\}$ for $i \in \mathcal{G}_d$: commitment status of discretionary generator $i$. $u_i=0$ is off and $u_i=1$ is on.
    \item $v_i \in \{0, 1\}$ for $i \in \mathcal{G}_d$ when $u_i$ was previously 1. $v_i=0$ represents staying off and $v_i=1$ is starting to ramp up.
\end{itemize}

The parameters supplied to the optimization solver are:
\begin{itemize}
    \item $a_i, b_i, c_i$: quadratic cost coefficients usually resulting from fuel costs or strategic bidding behavior.
    \item $P_{\max,i}, P_{\min,i}$: maximum capacity and minimum stable output, determined by a unit's physical characteristics.
    \item $D$, $\sigma_D$: demand volume and the possible demand deviations in the robust model formulation.
 \item $K_i$: commitment change penalty for generator $i$. 
    \item $r_{d,i}$ and $r_{u,i}$: generator $i$'s  limits for ramp down normalized so that $\Delta t = 1$.
    \item $\beta$: bias factor for selecting more efficient generators to meet the contingencies. In our model we set it to $0.001$. 
    \item $P_{typ,i}$, a typical production level of a generator during a peak event, which we define as $(4P_{\max,i}+P_{\min,i})/5$.
    \item $P_{t-1,i}$ - production level for the generator at the previous time period. All variables are assumed to refer to the current time, unless otherwise noted.
    \item $S_r$ - supply produced by generators that are already ramping before the solve (they cannot be changed)\footnote{For generators ramping up, they are assumed to generate nothing for some number of periods, and then increase linearly to $P_{\min,i}$ in some number of time periods specific to the generator. Generators ramping down are assumed to decrease linearly from $P_{\min,i}$ in some number of time periods specific to the generator.}.
    \item $S_{\max,r}$ total capacity of ramping up units  (ramping down generators will be off before the contingency)
    \item $S_{\min,r}$ minimum working capacity of ramping up units 
    \item ${D}_{\min,72}$ and ${D}_{\max,72}$: minimum and maximum volumes of predicted demand in next 72 hours.
    \item $\mathcal{G}_m$: the set of must run generators.
    \item $\mathcal{G}_d$: the set of discretionary generators.
    \item $R = \max\{\max_{i \in \mathcal{G}_m} P_{\max,i}, \max_{i \in \mathcal{G}_d} P_{\max,i}\}$: reserve margin (changed to only include generators that can run).
\end{itemize}

\cref{alg:ramp} summarizes the workflow for this section.
\begin{algorithm}[htbp]
\caption{Runtime and ramp constrained RRUC.}
\label{alg:ramp}
\begin{algorithmic}[1]
\For{$i=1$ \textbf{to} number\_of\_periods}
\State Update rolling 24 hour  generator states (use them in steps 3-8).
  \State Increment ramping periods for ramping units

\State Switch unit that finished preparation to ramping up
  \State Turn units that have finished ramping up to on
  \State {Turn units that have finished ramping down to off}
    \State {Find must run units (with minimum on time not yet reached or with $P_{t-1,i}-P_{min,i}$ too large to shut off)}
  \State Find can't start generators (those that have exceeded  maximum allowed daily starts or are already ramping)
\State Solve \cref{eq:OptimizationEq5} with $u_i,v_i\in\{0,1\}\rightarrow y_i,z_i\in[0,1]$
  \State  Order units in $\mathcal{G}_d$ into sets descending by $y_i$, then $z_i$
\For{Each $\mathcal{G}_d$ satisfying (\cref{eq:cont5})-(\cref{eq:MinConst5})} 
\State Solve economic dispatch \cref{eq:OptimizationEco} with $\mathcal{G}_m, \mathcal{G}_d$
\If{Objective is lowest so far}
\State Update best solution and objective
\EndIf
\EndFor
\State \textbf{save} Best solution and best objective
\EndFor
\end{algorithmic}
\end{algorithm}

\subsection{Simulation and Results}

Much of the data used in this Section has already been presented in \cref{sec:runtime} for runtime constraints. The additional assumptions about equipment characteristics are presented in \cref{tab:dataramp}. While ramp characteristics usually come in MW/min, we use the \% of $P_{\max,i}$ per minute as the original input into the simulation. We compile the percentage values the from the academic reports of Deutsches Institut für Wirtschaftsforschung \cite{diw2013} and the manufacturer catalog of GE Vernova~\cite{GE}. The provided numbers, although summarising different sources, may still present the upper bound of ramp rates, as comprehensive survey of large thermal units suggest that in real systems power plants do not always deliver the nominal ramp rates specified for equipment~\cite{India}.

Not all the available generating units could be accurately matched to the PJM set of generators on $P_{\min,i}$ and $P_{\max,i}$. Therefore we used approximations based on size and the assumed fuel used by a turbine (gas or coal steam). Larger turbines were assigned lower ramp values, as were coal steam turbines. Once the ramp percentages were assigned to specific generating units, the $P_{\max,i}$ of those generators was multiplied by the ramp rate to get the MW/min value, which was later used in the simulation.

If the system starts with all generators off, \cref{eq:OptimizationEq5} is usually infeasible. We initialized the system with the fastest ramping generators that can supply the previous demand. This may be suboptimal, but it is not the focus of this paper. Operators can initialize the system based on their own cold start protocol to a state that is known to be stable. If operators are warm starting using \cref{alg:ramp} they can use the current system state.

\begin{table}[htb]
\caption{\label{tab:dataramp} Ramp Assumptions and Data Sources}
\begin{tabular}{lrrr} \toprule
   Parameter & Source & Coal & Gas
\\ \midrule
Ramp up rate [\%/min]  & ~\cite{diw2013}, ~\cite{GE} & 1-6 & 2-12
\\ Ramp down rate [\%/min] & ~\cite{diw2013} & 5 & 15
\\ \bottomrule
\end{tabular}
\end{table}

\begin{figure}[htbp]
\centering
\includegraphics[width=0.95\linewidth]{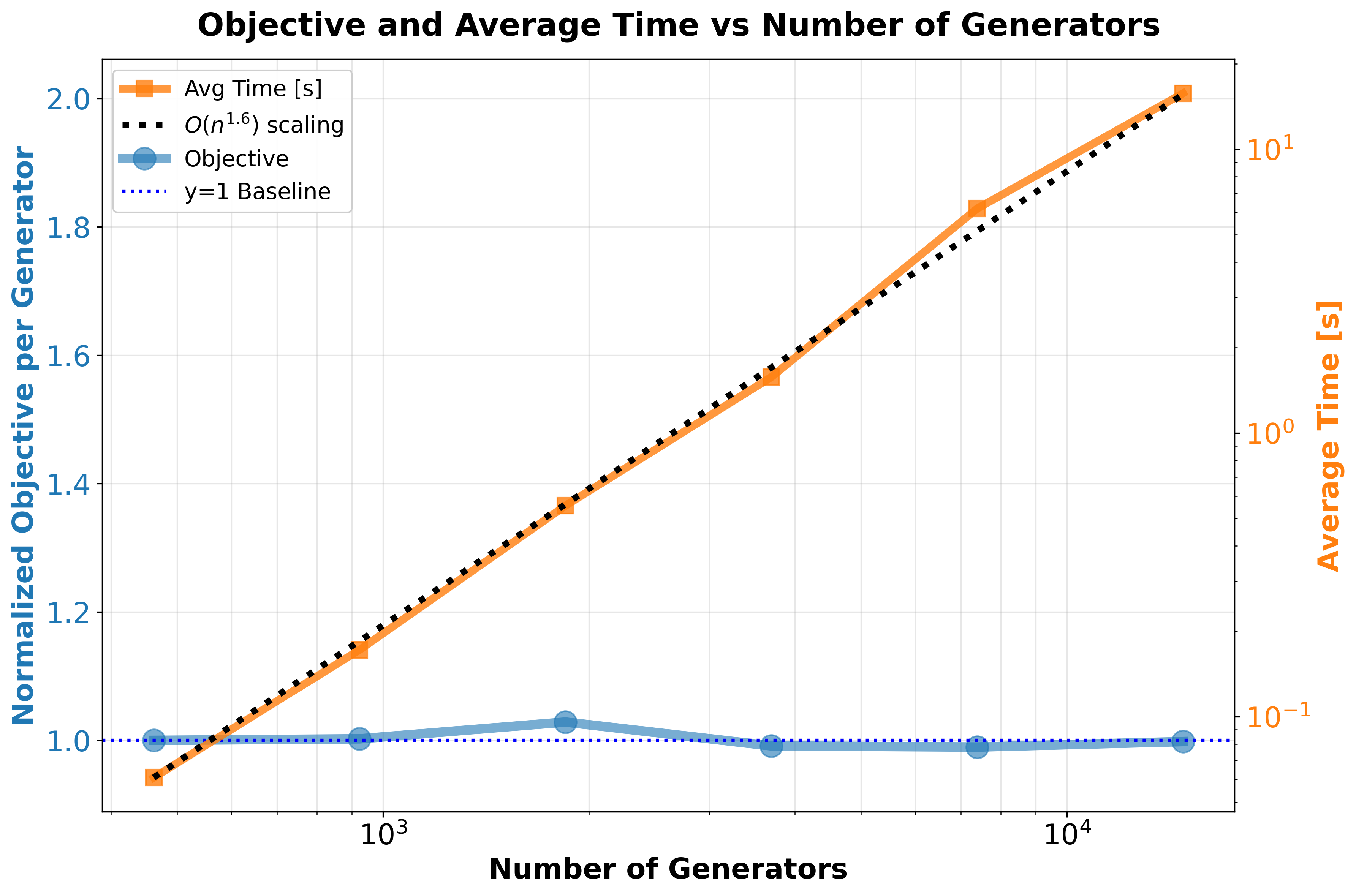}
\caption{\label{fig:ramp_time_obj} Runtime and objective function value for RRUC algorithm with runtime constraints and ramping constraints for 462-14784 generators with the \cref{fig:ramp_profile} profile. The objective function per generator fluctuates, but does not show an increasing or decreasing trend based on the system size. The runtime scales as $\approx O(n^{1.6})$, where $n$ is the number of generators. }
\end{figure}

\cref{fig:ramp_time_obj} summarizes the runtime and objective results. These results are comparable to those found in \cref{sec:runtime} using the same software and hardware.  The runtime scaled roughly as $O(n^{1.6})$, where $n$ is the number of generators. Ideal scaling is believed to be bounded below by $O(n^{1.5})$, due to the scaling of the underlying linear system solvers~\cite{morrison2016}. This suggests that our optimization solver has close to ideal scaling. It is well below the exponential cost of MIP solvers. The average cost of the system did not increase with the system size. So adding extra generators does not result in scale cost effects.

\section{Unit Commitment With Smooth Ramp Up}
\label{sec:nonconst}

Finally, we tested the suggested RRUC methodology on systems with smooth ramp up. The use of constant ramp rates is common in research, but they represent a simplified view of generator constraints. The ramp rates are not identical for ramping up and ramping down\cite{diw2013}. Ramping rates are also not constant during the ramping process, later stages of ramping are different from earlier stages\cite{chen2023}, \cite{Song2014}. Further, there is no clear distinction in practice between the prepare and ramp up stages in \cref{fig:ramp_profile}. They are both part of a deterministic shift from the generator being off to on. 

\begin{figure}[htbp]
\centering
\includegraphics[width=0.95\linewidth]{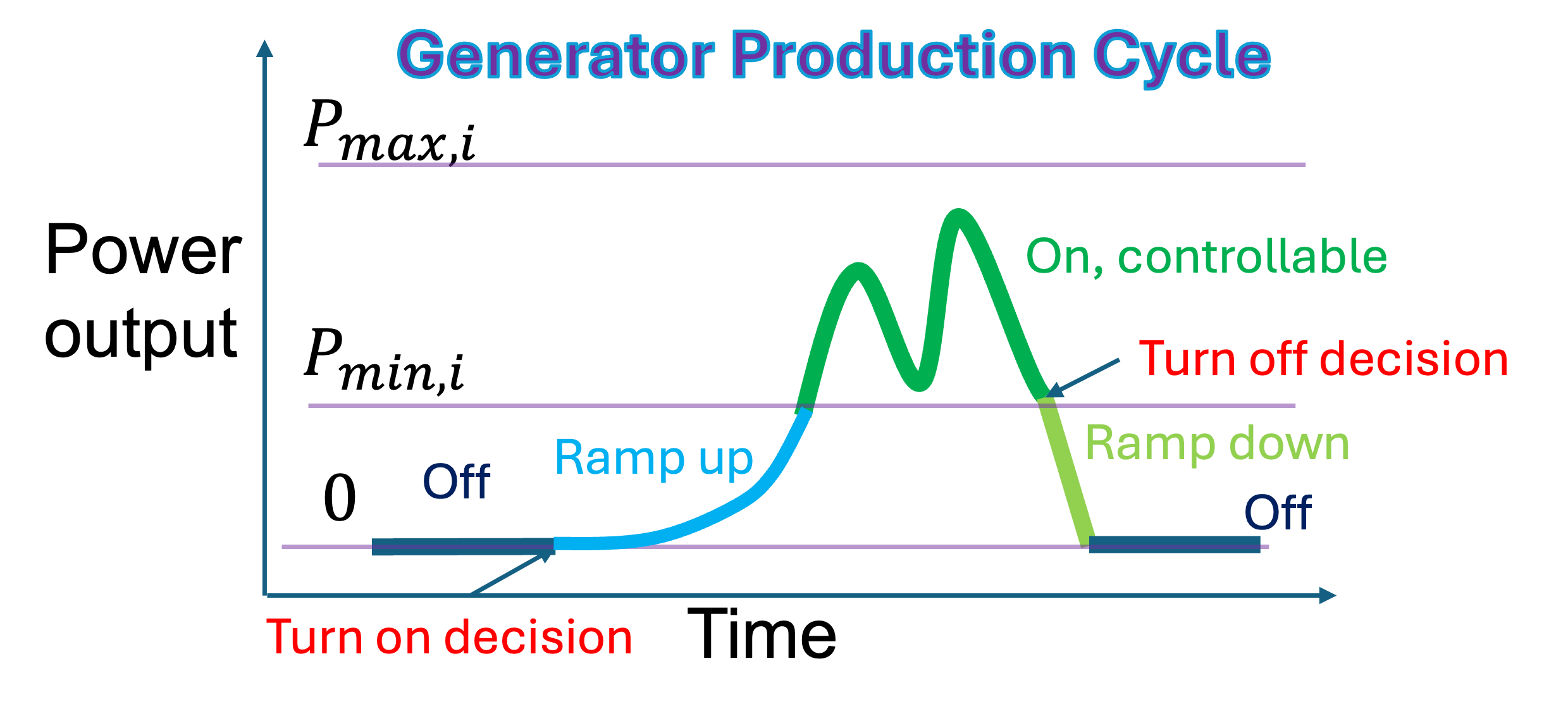}
\caption{\label{fig:ramp_profile_smooth} Unit $i$'s typical production cycle in \cref{sec:nonconst}. When on, $P_i\in [P_{\min,i},P_{\max,i}]$, though it cannot exceed its ramping power. It can turn off if $P_{t-1,i}-r_{d,i}\leq P_{\min,i}$. The only other decision can be made when the unit is off, to start ramping up. During other periods, it is uncontrollable.}
\end{figure}

We introduce a model extension in \cref{fig:ramp_profile_smooth} that accounts for these discrepancies to create a more realistic model of generator states. We combine the prepare and ramp up stages into one stage where production increases quadratically from $0$ to $P_{\min,i}$. This produces a more realistic, smoothed version of \cref{fig:ramp_profile}. Due to rounding up to whole time periods, this has the added effect of potentially requiring one less period for ramping up, which is closer to the actual ramp up time.

\begin{figure}[htbp]
\centering
\includegraphics[width=0.95\linewidth]{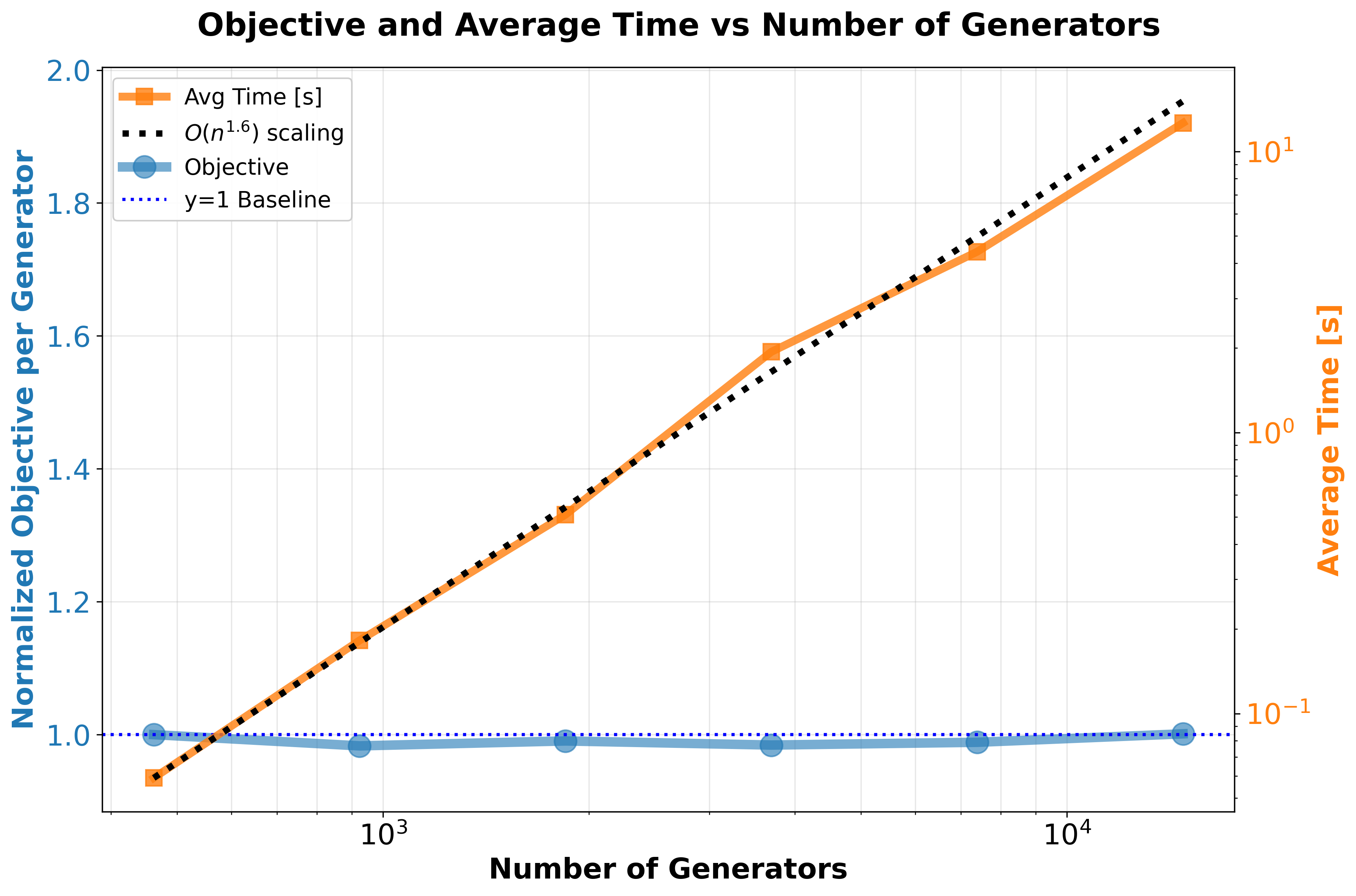}
\caption{\label{fig:ramp_time_obj_smooth} Runtime and objective function value for RRUC algorithm with runtime constraints and ramping constraints for 462-14784 generators, with the \cref{fig:ramp_profile_smooth} profile. The objective function per generator fluctuates, but does not show an increasing or decreasing trend based on the system size. The runtime scales as $\approx O(n^{1.6})$.}
\end{figure}

\cref{fig:ramp_time_obj_smooth} shows results using the \cref{fig:ramp_profile_smooth} generator profile model. The results were similar to those of \cref{sec:ramp}. There is no upward or downward trend in objective value per generator as the number of generators grows. Runtime continued to scale sub-quadratically, slightly slower than the optimal $O(n^{1.5})$. 
\begin{figure}[htbp]
\centering
\includegraphics[width=0.95\linewidth]{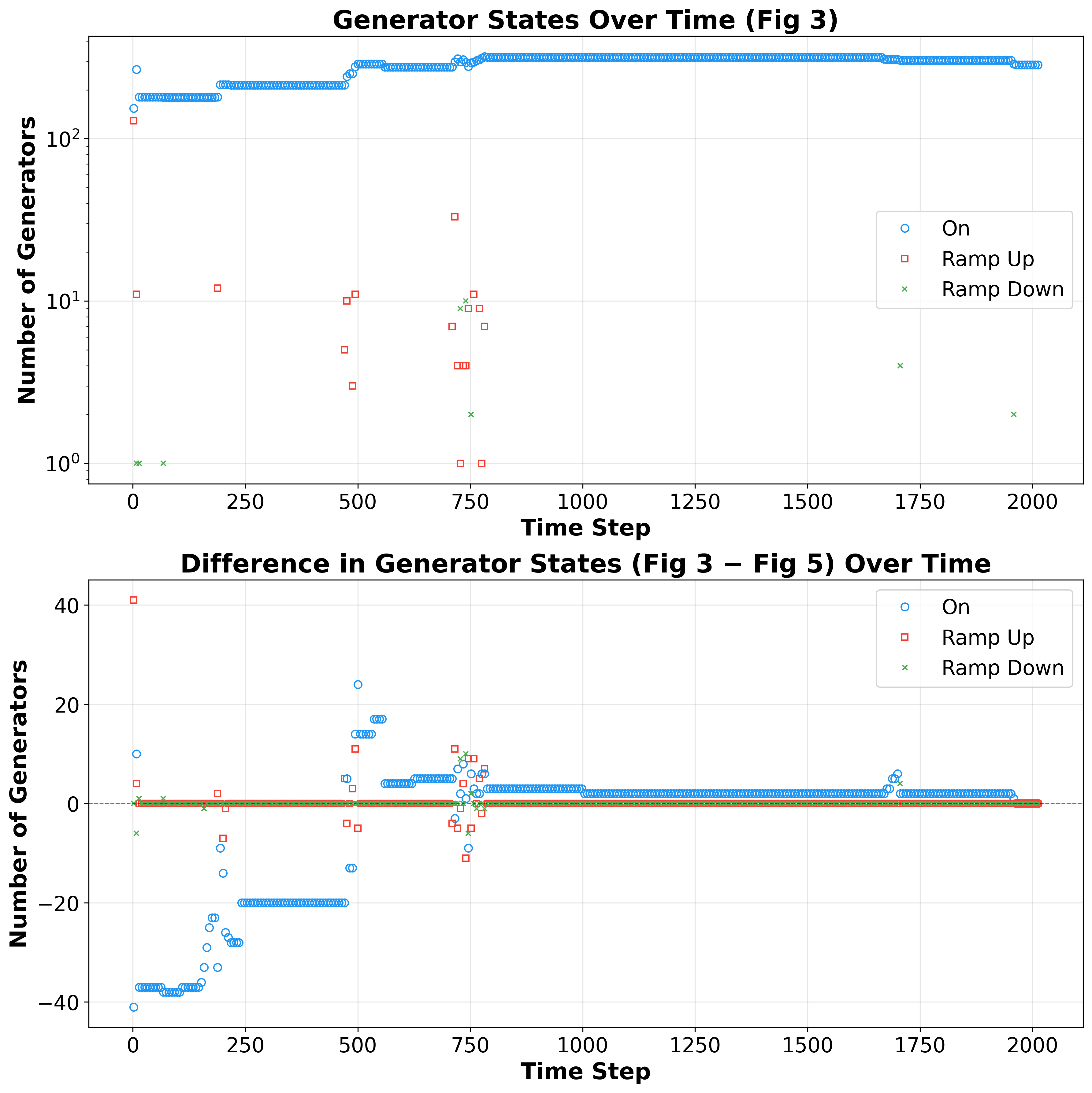}
\caption{\label{fig:gen_state} Number of generators in a given state for the 462 generator simulation. Ramp up includes the preparation stage for \cref{sec:ramp} results. Note: Only every $6$th time step is shown for readability. (Top) Total number of generators in the \cref{sec:ramp} simulation in on, ramp up, and ramp down states. Where there are $0$ generators in a given state, it is omitted on the log-scale. At almost all time steps, the majority of active generators are on and controllable.  (Bottom) Difference in the number of generators in a given state between the simulation in \cref{sec:ramp} (\cref{fig:ramp_profile} profile) and \cref{sec:nonconst} (\cref{fig:ramp_profile_smooth} profile). At almost all time steps, there is no difference.}
\end{figure}

The top part of \cref{fig:gen_state} shows that most time steps have no generators ramping up or down. In $1881$ out of $2016$ time steps ($93.3\%$), all generators were on or off for the \cref{fig:ramp_profile} generator model. In $1884$ out of $2016$ time steps ($93.5\%$), all generators were on or off for the \cref{fig:ramp_profile_smooth} generator model. Ramp up and ramp down bursts tend to be short and include many generators. The bottom part of \cref{fig:gen_state} shows the discrepancies between the solutions that both generator profiles produce. The differences in the generator states of the solutions remain large for the first $\approx 800$ time steps (equivalent to almost $3$ days), but level off. Small differences in the input can lead to very different but similar quality (similar cost) solutions, when aggregated over all time periods. This shows the non-convexity of the problem and the strength of the optimization solver, which is able to find different minima.

\section{Conclusion and Further Work}

We proposed an algorithm to solve unit commitment with mechanical and ramping constraints by relaxing the integer constraints and then rounding solutions to satisfy all constraints. Our algorithm shows $\approx O(n^{1.6})$ time scaling in the number of generators. As the problem size increases, it retains accuracy nearly perfectly.

For a system with only runtime constraints, a $2\times$ increase in the number of units results in about a $2.6\times$ increase in the computational time on average ($\approx O(n^{1.37})$ scaling). Further, the average value of the objective function increased by less than $1.05\times$ for each $2\times$ increase in the system size, indicating that there are some small diseconomies of scale associated with the increase in system size using the suggested algorithm. 

When ramping constraints were added to the problem, runtime scaling was close to ideal. For every $2\times$ increase in problem size, runtime increased about $3\times$ on average.  The  objective (cost) function per generator did not increase with the problem size. It is possible that being more restrictive is beneficial for a ramp constrained formulation for objective scaling. The demand at each time step is close to the one preceding it, so a more restricted solver with ramping constraints can restrict the search on better solutions. This means that even as the problem grows, at the cost of more runtime, the solver can pinpoint similar quality solutions.

A switch to a smoother generator profile did not change the average results drastically. Objectives and their perfect scaling were similar. Runtimes and their scaling were slightly better. For every $2\times$ increase in the problem size, the runtime increased less than $3\times$ on average. However, depending on the generator profile, certain time periods have substantially different units in different states. This indicates that even slightly different ramping profiles can lead to very different unit commitments, even if their objectives are similar. This highlights the non-convexity of the problem and the importance of accurate knowledge about generator profiles.

Further research would be required to assess the performance of the proposed methods when applied to a fully temporally coupled objective function and constraints. Further, it would be valuable to identify if the suggested algorithm demonstrates good convergence and scaling in the presence of transmission constraints.

\label{sec:summary}

\section*{Acknowledgements}
The authors thank Nicholson Koukpaizan for contributing to an improved manuscript with his thorough technical review.
\bibliographystyle{IEEEtran}
\bibliography{bibfile,refs-slaven}

\end{document}